\documentclass[12pt,a4paper]{article}
\usepackage[cp1251]{inputenc} 
\usepackage[russian]{babel}
\usepackage{amsfonts}

\newcommand{\di}{\displaystyle}

\newcommand{\al}{\alpha}
\newcommand{\be}{\beta}
\newcommand{\ga}{\gamma}
\newcommand{\de}{\delta}
\newcommand{\la}{\lambda}

\newcommand{\vv}{\varphi}
\newcommand{\iy}{\infty}

\begin{document}

\begin{center}
{\large\bf
Recovering First Order Integro-Differential Operators from Spectral Data}\\[0.3cm]
{\bf V.A.\,Yurko} \\[0.5cm]
\end{center}

\thispagestyle{empty}

\noindent {\bf Abstract.} First order integro-differential operators
on a finite interval are studied. Properties of spectral characteristic
are established, and the uniqueness theorem is proved for the
inverse problem of recovering operators from their spectral data.



{\bf 1. } Inverse spectral problems consist in recovering operators
from their spectral characteristics. Such problems often appear in
mathematics, mechanics, physics, electronics, geophysics and
other branches of natural sciences and engineering. The greatest
success in the inverse problem theory has been achieved for the
Sturm-Liouville operator (see, e.g., [1-3]) and afterwards for
higher-order differential operators [4-6] and other classes
of differential operators.

For integro-differential and other classes of nonlocal operators inverse
problems are more difficult for investigation, and the main classical methods
(transformation operator method and the method of spectral mappings [1-6])
either are not applicable to them or require essential modifications, and
for such operators the general inverse problem theory does not exist.
At the same time, nonlocal and, in particular, integro-differential operators
are of great interest, because they have many applications (se, e.g., [7]).
We note that some aspects of inverse problems for integro-differential
operators were studied in [8-10] and other works. In the present paper
we study inverse spectral problem for one class of first order
integro-differential operators on a finite interval. Properties of
spectral characteristic are established, and the uniqueness theorem
is proved for this class of inverse problems.

\bigskip
{\bf 2.} Consider the integro-differential equation
$$
\ell y:=iy'(x)+\int_0^x Q(x,t)y(t)\,dt=\la y(x),\quad x\in[0,\pi],  \eqno(1)
$$
where $Q(x,t)=R(x)V(t).$ We assume that $R(x), V(t)$ are continuous
complex-valued functions, and
$$
R(\pi-x)\sim C_\al x^\al,\; V(x)\sim D_\be x^\be,\;
x\to+0,\quad C_\al D_\be\ne 0.
$$
Let $\vv(x,\la)$ be the solution of Eq. (1) with the condition
$\vv(0,\la)=1.$ Then the following representation holds (see [3]):
$$
\vv(x,\la)=\exp(-i\la x)+\int_0^x K(x,t)\exp(-i\la t)\,dt,            \eqno(2)
$$
where $K(x,t)$ is a continuous function, and $K(x,0)=0.$ Denote
$$
\Pi_{+}:=\{\la:\; \mbox{Im}\,\la\ge 0\},\quad
\Pi_{-}^{\de}:=\{\la:\; \arg\la\in [\pi+\de,2\pi-\de]\}.
$$
It follows from (2) that for $|\la|\to\iy$ uniformly in $x\in[0,\pi]$:
$$
\left.\begin{array}{c}
\vv^{(\nu)}(x,\la)=(-i\la)^{\nu}\exp(-i\la x)(1+o(1)),
\quad \la\in\Pi_{+},\quad \nu=0,1,\\[3mm]
\vv^{(\nu)}(x,\la)=o(\la^\nu),\quad \la\in\Pi_{-}^{\de},\quad\nu=0,1.
\end{array}\right\}                                                 \eqno(3)
$$
Denote
$$
\vv_{\nu}(x,\la):=\di\frac{1}{\nu!}
\di\frac{\partial\vv(x,\la)}{\partial\la^\nu},\;\nu\ge 0,\quad
\Delta(\la):=\vv(\pi,\la).
$$
The function $\Delta(\la)$ is entire in
$\la$ of exponential type, and its zeros $\Lambda:=\{\la_n\}_{n\ge 1}$
(counting with multiplicities) coincide with the eigenvalues of the
boundary value problem $L=L(R,V)$ for Eq. (1) with the condition
$y(\pi)=0.$ Let $m_n$ be the multiplicity of $\la_n$
($\la_n=\la_{n+1}=\ldots = \la_{n+m_n-1}$). Denote
$$
S:=\{n:\; n-1\in{\bf N},\; \la_{n-1}\ne\la_{n}\}\cup\{1\},\quad
s_{n+\nu}(x):=\vv_{\nu}(x,\la_n),\; n\in S,\; \nu=\overline{m_n-1}.
$$
The functions $\{s_n(x)\}_{n\ge 1}$ are eigen and associated functions
for $L.$

\smallskip
{\bf Example 1. } Let $\la_1=\la_2<\la_3<\la_4=\la_5=\la_6<\la_7<
\la_8<\ldots$. Then $S=\{1,3,4,7,8,\ldots\},$ $s_1(x)=\vv_0(x,\la_1),$
$s_2(x)=\vv_1(x,\la_1),$ $s_3(x)=\vv_0(x,\la_3),$ $s_4(x)=\vv_0(x,\la_4),$ $s_5(x)=\vv_1(x,\la_4),$ $s_6(x)=\vv_2(x,\la_4),$
$s_7(x)=\vv_0(x,\la_7),\ldots$

\smallskip
Let the function $\eta(x,\la)$ be the solution of the problem
$$
i\eta'(x,\la)-R(x)\int_x^\pi V(t)\eta(t,\la)\,dt
+R(x)=\la\eta(x,\la), \quad \eta(\pi,\la)=0.                        \eqno(4)
$$
Denote $\theta(x,\la):=\eta(\pi-x,\la).$ Then
$$
i\theta'(x,\la)+R_0(x)\int_0^x V_0(t)\theta(t,\la)\,dt
-R_0(x)=-\la\theta(x,\la), \quad \theta(0,\la)=0,                   \eqno(5)
$$
where $R_0(x):=R(\pi-x),\,V_0(x):=V(\pi-x).$ It follows from (5) that
$$
\theta(x,\la)=\int_0^x g(x,t,\la)R_0(t)\,dt,                        \eqno(6)
$$
where $g(x,t,\la)$ is  Green's function of the Cauchy problem, and
$$
ig_x(x+t,t,\la)-\la g(x+t,t,\la)+R(x+t)
\int_0^x V(\tau+t)g(\tau+t,t,\la)\,d\tau=0,\; g(t,t,\la)=-i,
$$
and consequently, $g(x+t,t,\la)=-i\vv(x,\la;t),$ where
$\vv(x,\la;t)$ is the solution of the Cauchy problem
$$
i\vv'(x,\la;t)+R(x+t)\int_0^{x} V(\tau+t)\vv(\tau,\la;t)\,d\tau
=\la \vv(x,\la;t),\; \vv(0,\la;t)=1.
$$
In view of (2) we get
$$
\vv(x,\la;t)=\exp(-i\la x)+\int_0^x K(x,\tau;t)\exp(-i\la\tau)\,d\tau,
$$
where $K(x,\tau;t)$ is a continuous function. This yields
$$
g(x,t,\la)=-i\exp(-i\la(x-t))
-i\int_0^{x-t} K(x-t,\tau;t)\exp(-i\la\tau)\,d\tau.                \eqno(7)
$$
Substituting (7) into (6), we obtain
$$
\theta(x,\la)=\int_0^x P(x,t)\exp(i\la t)\,dt,                     \eqno(8)
$$
where
$$
P(x,t)=-iR_0(x-t)-i\int_0^{x-t} R_0(\tau)K(x-\tau,t;\tau)\,d\tau.  \eqno(9)
$$
Clearly, $P(x,x)=-iR_0(0),\; P(x,0)=-iR_0(x).$ Using (8)-(9) and (4)
we conclude that for $|\la|\to\iy$ uniformly in $x\in[0,\pi]$:
$$
\left.\begin{array}{c}
\eta^{(\nu)}(x,\la)=o(\la^{\nu}),
\quad \la\in\Pi_{+},\quad \nu=0,1,\\[3mm]
\eta^{(\nu)}(x,\la)=o(\la^\nu\exp(i\la(\pi-x))),
\quad \la\in\Pi_{-}^{\de},\quad\nu=0,1.
\end{array}\right\}                                               \eqno(10)
$$
Denote
$$
\Delta_0(\la):=1-\int_0^\pi V(t)\eta(t,\la)\,dt.                   \eqno(11)
$$
Using (4) and (11) we calculate
$$
i\eta'(x,\la)+R(x)\Big(\Delta_0(\la)+\int_0^x V(t)\eta(t,\la)\,dt\Big)
=\la\eta(x,\la), \quad \eta(\pi,\la)=0.                           \eqno(12)
$$
In particular, it follows from (12) that zeros of the entire
function $\Delta_0(\la)$ coincide with the zeros of $\Delta(\la),$
and multiplicities of zeros of $\Delta_0(\la)$ are not more than
multiplicities of zeros of $\Delta(\la).$ Therefore the function
$\Delta(\la)/\Delta_0(\la)$ is entire in $\la$ of exponential type.
Denote $\Delta_1(\la):=\Delta_0(\la)\exp(-i\la\pi).$ Using (8), (9) and
(11), by standard arguments (see, for example, [?]) we obtain that
for $|\la|\to\iy,$ the following asymptotical formulae hold
$$
\left.\begin{array}{c}
\Delta_1(\la)=\exp(-i\la\pi)(1+o(1)),\quad \la\in\Pi_{+},\\[3mm]
\Delta_1(\la)=B\la^{-\ga-1}(1+o(1)),\quad \la\in\Pi_{-}^{\de},
\end{array}\right\}                                               \eqno(13)
$$
where $B\ne 0,\; \ga:=\al+\be+1.$
The function $F(\la):=\Delta(\la)/\Delta_1(\la)$ is entire in
$\la$ of exponential type. By virtue of (3),
$$
\left.\begin{array}{c}
\Delta(\la)=\exp(-i\la\pi)(1+o(1)),\quad \la\in\Pi_{+},\\[3mm]
\Delta(\la)=o(1),\quad \la\in\Pi_{-}^{\de}.
\end{array}\right\}                                              \eqno(14)
$$
Together with (13) this yields that $F(\la)\equiv 1,$ i.e.
$\Delta(\la)\equiv\Delta_1(\la)$ or
$$
\Delta(\la)\equiv\Delta_0(\la)\exp(-i\la\pi).                   \eqno(15)
$$
Denote
$$
\eta_{\nu}(x,\la):=\di\frac{1}{\nu!}
\di\frac{\partial\eta(x,\la)}{\partial\la^\nu},\;\nu\ge 0,\;
\psi_{n+\nu}(x):=\eta_{\nu}(x,\la_n),\; n\in S,\;\nu=\overline{m_n-1}.
$$
The functions $\{\psi_n(x)\}_{n\ge 1}$ are eigen and associated
functions for the boundary value problem $L,$ and
$$
\psi_{n+\nu}(x)=\sum_{j=0}^\nu \be_{n+\nu-j}s_{n+j}(x),
\; n\in S,\;\nu=\overline{m_n-1}.                             \eqno(16)
$$
The coefficients $\{\be_n\}_{n\ge 1}$ are called Levinson's
weight numbers, and the data $\{\la_n,\be_n\}_{n\ge 1}$ are
called the spectral data for the boundary value problem $L.$
We will consider the following inverse problem:

\smallskip
{\bf Inverse problem 1. }{\it Given the spectral
data $\{\la_n,\be_n\}_{n\ge 1}$, construct $R$ and $V.$}

\medskip
{\bf 3.} Below we will assume that $R(x)\ne 0$ a.e. on $(0,\pi).$
If this condition does not hold, then the specification of the
spectral data does not uniquely determine $L$ (see Example 2 below).

Let us formulate the uniqueness theorem for this inverse problem.
For this purpose, together with $L$ we consider the boundary
value problem $\tilde L:=L(\tilde R, \tilde V)$ of the same
form but with a different functions $\tilde R(x), \tilde V(t).$
We agree that everywhere below if a certain symbol $\al$ denotes
an object related to $L,$ then $\tilde\al$ will denote the
analogous object related to $\tilde L.$

\smallskip
{\bf Theorem 1. }{\it Let $\{\tilde\la_n,\tilde\be_n\}$ be the
spectral data for the problem $\tilde L=L(\tilde R,\tilde V).$
If $\la_n=\tilde\la_n$, $\be_n=\tilde\be_n$ for all $n\ge 1,$
then $R(x)\equiv\tilde R(x),$ $V(x)\equiv\tilde V(x),$
$x\in[0,\pi].$}

\smallskip
{\it Proof. } Using (14)-(15) and Hadamard's factorization
theorem we get $\Delta_0(\la)\equiv\tilde\Delta_0(\la).$
Taking (16) into account, we deduce that the functions
$$
A_j(x,\la)=(\Delta_0(\la))^{-1}\exp(i\la x)
\Big(\tilde\vv(x,\la)\eta^{(j-1)}(x,\la)
-\tilde\eta(x,\la)\vv^{(j-1)}(x,\la)\Big),\; j=1,2,
$$
are entire in $\la$ of exponential type. Taking (3), (10) and
(13) into account we obtain for $|\la|\to\iy$:
$$
A_1(x,\la)=o(1),\; A_2(x,\la)=o(\la),\quad \la\in\Pi_{+},
$$
$$
A_1(x,\la)=o(\la^{\ga+1}),\; A_2(x,\la)=o(\la^{\ga+2}),
\quad \la\in\Pi_{-}^{\de},
$$
and consequently,
$$
A_1(x,\la)\equiv 0,\quad A_2(x,\la)\equiv A(x),                \eqno(17)
$$
where the function $A(x)$ does not depend on $\la.$
In particular, (17) yields
$$
\tilde\vv(x,\la)\eta(x,\la)\equiv\tilde\eta(x,\la)\vv(x,\la), \eqno(18)
$$
$$
\tilde\vv(x,\la)\eta'(x,\la)-\tilde\eta(x,\la)\vv'(x,\la)
\equiv A(x)\Delta_0(\la)\exp(-i\la x).                        \eqno(19)
$$
Similarly, we obtain
$$
\vv(x,\la)\eta'(x,\la)-\eta(x,\la)\vv'(x,\la)
\equiv A^{*}(x)\Delta_0(\la)\exp(-i\la x),                    \eqno(20)
$$
where $A^{*}(x)$ does not depend on $\la.$ Using (18) we calculate
$$
\tilde\vv(x,\la)\Big(\vv(x,\la)\eta'(x,\la)-\eta(x,\la)\vv'(x,\la)\Big)=
\vv(x,\la)\Big(\tilde\vv(x,\la)\eta'(x,\la)-\tilde\eta(x,\la)\vv'(x,\la)\Big).
$$
Together with (19)-(20) this yields
$$
\tilde\vv(x,\la)A^{*}(x)\equiv\vv(x,\la)A(x).
$$
Taking (3) into account, we conclude that $A(x)\equiv A^{*}(x),$ and
$$
(\tilde\vv(x,\la)-\vv(x,\la))A(x)\equiv 0.                    \eqno(21)
$$
Furthermore, using (18), (19) and equations (1) and (4), we infer
$$
iA(x)\Delta_0(\la)\exp(-i\la x)=i\tilde\vv(x,\la)\eta'(x,\la)
-i\tilde\eta(x,\la)\vv'(x,\la)
$$
$$
=\tilde\vv(x,\la)\Big(R(x)\int_x^\pi V(t)\eta(t,\la)\,dt-R(x)\Big)
+\tilde\eta(x,\la)R(x)\int_0^x V(t)\vv(t,\la)\,dt.
$$
Hence, for $|\la|\to\iy,\; \la\in\Pi_{+}$, we get $A(x)\equiv iR(x).$
In view of (21), one has
$$
(\tilde\vv(x,\la)-\vv(x,\la))R(x)\equiv 0.                   \eqno(22)
$$
Since $R(x)\ne 0$ a.e. on $(0,\pi),$ it follows from (22) that
$\tilde\vv(x,\la)\equiv\vv(x,\la).$ By virtue of (18),
$\tilde\eta(x,\la)\equiv\eta(x,\la).$ Then, according to (4),
$$
R(x)-\tilde R(x)=R(x)\int_x^\pi V(t)\eta(t,\la)\,dt-
\tilde R(x)\int_x^\pi \tilde V(t)\eta(t,\la)\,dt.
$$
For $|\la|\to\iy$ this yields $R(x)\equiv\tilde R(x),$ $x\in[0,\pi],$
and consequently $V(x)\equiv\tilde V(x),$ $x\in[0,\pi].$
Theorem 1 is proved.

\smallskip
{\bf Example 2. } Fix $a\in(0,\pi).$ Let $R(x)\equiv 0$ for $x\in [0,a],$
and $R(x)\ne 0$ for $x\in(a,\pi).$ Put $\tilde R(x)\equiv R(x)$ for
$x\in[0,\pi],$ and chose $V(t), \tilde V(t)$ such that
$V(t)\equiv\tilde V(t)$ for $t\in (a,\pi),$ and
$V(t)\ne\tilde V(t)$ for $t\in [0,a].$ Then
$\tilde\vv(x,\la)\equiv\vv(x,\la)$ and $\tilde\eta(x,\la)\equiv\eta(x,\la)$;
hence $\tilde\la_n=\la_n$, $\tilde\be_n=\be_n$ for all $n\ge 1.$

\bigskip
{\bf Acknowledgment.} This work was supported by Grant 17-11-01193
of the Russian Science Foundation.

\begin{center}
{\bf REFERENCES}
\end{center}
\begin{enumerate}
\item[{[1]}] Marchenko V.A., Sturm-Liouville Operators and their Applications.
     Naukova Dumka,  Kiev, 1977;  English  transl., Birkh\"auser, 1986.
\item[{[2]}] Levitan B.M., Inverse Sturm-Liouville Problems. Nauka,
     Moscow, 1984; English transl., VNU Sci.Press, Utrecht, 1987.
\item[{[3]}] Freiling G. and Yurko V.A., Inverse Sturm-Liouville Problems
     and their Applications. NOVA Science Publishers, New York, 2001.
\item[{[4]}] Beals R., Deift P. and Tomei C.,  Direct and Inverse Scattering
     on the Line, Math. Surveys and Monographs, v.28. Amer. Math. Soc.
     Providence: RI, 1988.
\item[{[5]}] Yurko V.A. Method of Spectral Mappings in the Inverse Problem
     Theory. Inverse and Ill-posed Problems Series. VSP, Utrecht, 2002.
\item[{[6]}] Yurko V.A. Inverse Spectral Problems for Differential
     Operators and their Applications. Gordon and Breach, Amsterdam, 2000.
\item[{[7]}] Lakshmikantham V. and Rama Mohana Rao M. Theory of
     integro-differential equations. Stability and Control: Theory and
     Applications, vol.1, Gordon and Breach, Singapure, 1995.
\item[{[8]}] Yurko V.A., An inverse problem for integro-differential
     operators. Matem. zametki, 50, no.5 (1991), 134-146 (Russian);
     English transl. in Math. Notes, 50, no.5-6 (1991), 1188-1197.
\item[{[9]}] Kuryshova Yu. An inverse spectral problem for differential
     operators with integral delay. Tamkang J. Math. 42, no.3 (2011), 295-303.
\item[{[10]}] Buterin S.A. On the reconstruction of a convolution perturbation
     of the Sturm-Liouville operator from the spectrum, Diff. Uravn. 46 (2010),
     146--149 (Russian); English transl. in Diff. Eqns. 46 (2010), 150--154.
\end{enumerate}

\end{document}